\documentclass{article}
\usepackage{graphicx}

\usepackage{amsmath,amssymb,amsfonts,epsfig}

\newtheorem{prop}{Proposition}
\newtheorem{defi}[prop]{Definition}
\newtheorem{cor}[prop]{Corollary}
\newtheorem{rem}[prop]{Remark}

\newcommand{\proof}[1]{\vspace{1.5ex}\noindent{{\bf Proof:} #1 
}\vspace{1.5ex}}

\newcommand{\ex}{\medskip \noindent {\bf Example }}

\newcommand{\C}{\mathbb C}             
\newcommand{\R}{\mathbb R}             

\newcommand{\ts}[1]{\protect{{\textstyle{ #1}}}}
\newcommand{\mc}[1]{\protect{{\mathcal{#1}}}}
\newcommand{\id}{\protect{\rm{id}}}

\renewcommand{\deg}{{\rm ddeg }}

\title{The shuffle Hopf algebra and quasiplanar Wick products}

\author{Dorothea Bahns \thanks{Department Mathematik, Universit\"at Hamburg, Bundesstr. 55, D - 20146 Hamburg, Germany --- bahns@math.uni-hamburg.de}}

\date{October 15, 2007}

\begin{document}

\maketitle

\begin{abstract}
The operator valued distributions which arise in quantum field theory on the noncommutative Minkowski space can be symbolized by a generalization of chord diagrams, the dotted chord diagrams.  
In this framework, the combinatorial aspects of quasiplanar Wick products are understood in terms of the shuffle Hopf algebra of dotted chord diagrams, leading to an algebraic characterization of quasiplanar Wick products as a convolution. Moreover, it is shown that the distributions do not provide a weight system for universal knot invariants.
%
\end{abstract}


\section{Introduction}

\medskip
Tensor products of the operator valued distributions that appear in quantum field theory are in general ill-defined when pulled back to the diagonal, and the process of renormalization is necessary to define products which are well-defined distributions. For products of free field operators, the renormalization procedure leads to what is called the Wick product of quantum fields.

In~\cite{quasiI}, the quasiplanar Wick products were defined as a generalization of Wick products which is suitable for products of quantum fields on the noncommutative Minkowski space. The definition is based on a certain notion of locality and it is the first step towards a full renormalization theory in the Minkowskian noncommutative framework. This paper elaborates on the quasiplanar Wick products'  combinatorial and algebraic aspects, leaving the functional analytic aspects aside. Its first aim is to clarify that the graphs we used in~\cite{quasiI,bahnsDiss} to handle the combinatorics of quasiplanar Wick products, are a generalization of the classic chord diagrams studied e.g. in knot theory~\cite{kassel,barnatanIntro}. 

Chord diagrams carry a cocommutative Hopf structure, and it is natural to try to reformulate the combinatorial aspects of quasiplanar Wick products in this algebraic language, much in the spirit of~\cite{conneskreimer}. So, the paper's second aim is to reformulate the combinatorial aspects of quasiplanar Wick products as proved in~\cite{quasiI,bahnsDiss} in this algebraic setting and to show that the Hopf structure with the shuffle product and deconcatenation coproduct is the natural one in our context. A completely algebraic characterization of quasiplanar Wick products in terms of a convolution in the shuffle Hopf algebra is given in section~\ref{shuffleHopf}.
Section~\ref{qft} is devoted to explicitely relating these algebraic objects and relations to the operator valued distributions which arise in quantum field theory on the noncommutative Minkowski space. In the last section, it is shown that, although these distributions bear some similarity with weight systems~\cite{barnatanIntro}, they do not fulfill the 4T relation.


\section{Dotted chord diagrams}

\medskip 
Let us first recall the notion of a chord diagram. Let $L$ denote a directed simple polygonal arc in $[0,1]\times \R$, let $\partial L$ denote its boundary, that is, the set of its two endpoints. A chord diagram on $L$ is a finite set of ordered pairs of distinct points on $L\setminus \partial L$. The pairs of points are usually symbolized by connecting lines, the {\em chords} of the chord diagram, whose shape is irrelevant, e.g.
\begin{picture}(50,10)
\qbezier(20,0)(30,16)(40,0)
\put(20,0){\circle*{2}}
\put(40,0){\circle*{2}}
\put(10,0){\circle*{2}}
\put(30,0){\circle*{2}}
\qbezier(10,0)(20,16)(30,0)
\put(0,0){\vector(1,0){50}}
\end{picture}
for the set $\{(x_1,x_3), (x_2,x_4)\}$, where $x_1,x_2,x_3,x_4$ appear on the arc from left to right in that order.

{\defi \rm Let $L$ be an arc. A {\em dotted chord diagram} on $L$ is a finite set of ordered pairs of points on $L\setminus \partial L$.}

\medskip Observe that in this definition, we admit pairs $(x,x)$ of {\em nondistinct} points. We refer to such pairs as {\em dots}.
We continue to symbolize a pair of distinct points by its connecting line and symbolize a pair $(x,x)$  on the arc simply by the point $x$ itself, e.g. 
\begin{picture}(60,10)
\put(20,0){\circle*{2}}
\put(40,0){\circle*{2}}
\put(10,0){\circle*{2}}
\put(30,0){\circle*{2}}
\put(50,0){\circle*{2}}
\qbezier(10,0)(20,16)(30,0)
\qbezier(20,0)(35,16)(50,0)
\put(0,0){\vector(1,0){60}}
\end{picture}
for a set $\{(x_1,x_3), (x_2,x_5), (x_4,x_4)\}$ with $x_i\neq x_j$ for $i\neq j$ and where $x_1,x_2,x_3,x_4, x_5$ appear on the arc from left to right in that order.
Observe that a dotted chord diagram that only contains pairs of distinct points is indeed a chord diagram in the usual sense. These diagrams thus appear as special cases of dotted chord diagrams and will be called diagrams without dots. Likewise, we will call the dotted chord diagrams containing only pairs of non-distinct points diagrams without chords or chordless diagrams.

We say that a dotted chord diagram  has dot-degree $m$, and write $\deg(D)=m$, if its pairs are built from $m$ distinct points on the arc. For example, the two dotted diagrams above have dot-degree $4$ and $5$, respectively. Note that usually, the degree of a  chord diagram is the number of its chords, hence, the dot-degree of a chord diagram without dots is twice the ordinary degree.

Let $V_m$ denote the finite dimensional vector space that is spanned over $\C$ by all dotted chord diagrams of dot-degree $m$, and let $V=\bigoplus_{m\geq 0} V_m$ with $V_0=\C$.
 
It is well known that the vector space of chord diagrams forms a Hopf algebra, see for instance~\cite{barnatanIntro}. We will generalize this structure to the vector space of dotted chord diagrams. The unique way to glue together two directed arcs $L_1$ and $L_2$ (in that order from left to right), such that the resulting arc is again directed, extends by linearity to an associative product $\mu:V\otimes V \rightarrow V$. Usually, we write $a\cdot b$ or $ab$ for  $\mu(a\otimes b)$. The unit for this product is the empty diagram 
\begin{picture}(20,10)\put(5,3.5){\vector(1,0){10}}
\end{picture} which we also denote by $\emptyset$ or $1$. Now consider the coproduct $\Delta: V \rightarrow V \otimes V$, defined on diagrams as
\[
\Delta(D)=\sum_{\emptyset \subset D^\prime \subset D} D^\prime \otimes (D\setminus D^\prime)
\]
where the sum runs over all subdiagrams (including the empty diagram as well as $D$ itself). Its counit is the map $\epsilon$ that is equal to 1 on the empty diagram and 0 elsewhere, and the primitive elements in $V$ are 
\begin{picture}(15,10)
\put(5,0){\circle*{2}}
\put(0,0){\vector(1,0){15}}
\end{picture}
and 
\begin{picture}(25,10)
\put(5,0){\circle*{2}}
\put(15,0){\circle*{2}}
\qbezier(5,0)(10,16)(15,0)
\put(0,0){\vector(1,0){25}}
\end{picture}
It is standard to check that $(V,\mu, \Delta)$ is a bialgebra, and since $V$ is graded and connected, it follows that $V$ is a Hopf algebra. Its antipode is given inductively on the dot-degree of a diagram by $S(\emptyset)=\emptyset$, and for $D \neq \emptyset$,
\[
S(D)=-D-\sum_{\emptyset \subsetneq D^\prime \subsetneq D} S(D^\prime)\; \cdot\; (D\setminus D^\prime)
\]
We have, for example, 
\[
S(\begin{picture}(45,10)
\put(15,0){\circle*{2}}
\put(35,0){\circle*{2}}
\put(5,0){\circle*{2}}
\put(25,0){\circle*{2}}
\qbezier(5,0)(15,16)(25,0)
\put(0,0){\vector(1,0){45}}
\end{picture})
=
\begin{picture}(45,10)
\put(15,0){\circle*{2}}
\put(35,0){\circle*{2}}
\put(5,0){\circle*{2}}
\put(25,0){\circle*{2}}
\qbezier(15,0)(25,16)(35,0)
\put(0,0){\vector(1,0){45}}
\end{picture}
+
\begin{picture}(45,10)
\put(15,0){\circle*{2}}
\put(35,0){\circle*{2}}
\put(5,0){\circle*{2}}
\put(25,0){\circle*{2}}
\qbezier(15,0)(20,11)(25,0)
\put(0,0){\vector(1,0){45}}
\end{picture}
-
\begin{picture}(45,10)
\put(15,0){\circle*{2}}
\put(35,0){\circle*{2}}
\put(5,0){\circle*{2}}
\put(25,0){\circle*{2}}
\qbezier(25,0)(30,11)(35,0)
\put(0,0){\vector(1,0){45}}
\end{picture}
\]
$S$ is an algebra-antihomomorphism, that is, $S(ab)=S(b)S(a)$, and since $\Delta$ is cocommutative, we have $S^2=id_V$, see~\cite{kassel}.


\section{The Quasiplanar Wick map}

\medskip 
Let us recall and extend some definitions from~\cite{mmr} and \cite{quasiI}, respectively.
The labelled intersection graph of a chord diagram without dots is a graph whose vertices are the chords of $D$, numbered from $1$ to $\deg(D)/2$ in the order in which their starting points appear along the arc, and where two vertices are connected by an edge iff the corresponding two chords in $D$ intersect.

To extend this definition to dotted chord diagrams, we first establish how to label dotted chord diagrams. Let $D$ be a dotted chord diagram, then we label its chords and dots on the same footing  in the order as they appear along the arc, e.g. 
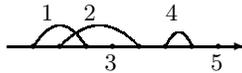
\[
\begin{picture}(90,25)
\put(10,7){\circle*{2}}
\put(20,7){\circle*{2}}
\put(30,7){\circle*{2}}
\put(40,7){\circle*{2}}
\put(37,-2){\small 3}
\put(50,7){\circle*{2}}
\put(60,7){\circle*{2}}
\put(70,7){\circle*{2}}
\put(80,7){\circle*{2}}
\put(77,-2){\small 5}
\qbezier(10,7)(20,23)(30,7)
\put(13,17){\small 1}
\qbezier(20,7)(35,23)(50,7)
\put(29,17){\small 2}
\qbezier(60,7)(65,18)(70,7)
\put(60,17){\small 4}
\put(0,7){\vector(1,0){90}}
\end{picture} 
\]

{\defi \rm The {\em labelled intersection graph} of a dotted chord diagram $D$  is a graph with coloured vertices. It is composed of the labelled intersection graph (with, say, white vertices) of the diagram $D$ with all dots removed but keeping the original labels of the diagram $D$, 
and an additional set of, say, black vertices, one for each dot in $D$, also with the labels from $D$. An edge connects a black vertex with a white vertex provided that on the arc, the dot corresponding to the black vertex is between the two endpoints of the chord corresponding to the white vertex. 
The adjacency matrix of the labelled intersection matrix of a dotted chord diagram was called the extended {\em incidence matrix} in~\cite{quasiI}.
A dotted chord diagram is called {\em connected} 
if its labelled intersection graph is connected. 
}

\ex The labelled intersection graph of the dotted chord diagram \linebreak
\begin{picture}(95,10)
\put(10,0){\circle*{2}}
\put(20,0){\circle*{2}}
\put(30,0){\circle*{2}}
\put(40,0){\circle*{2}}
\put(50,0){\circle*{2}}
\put(60,0){\circle*{2}}
\put(70,0){\circle*{2}}
\put(80,0){\circle*{2}}
\qbezier(10,0)(20,16)(30,0)
\qbezier(20,0)(35,16)(50,0)
\qbezier(60,0)(65,11)(70,0)
\put(3,0){\vector(1,0){87}}
\end{picture} 
is the graph
\begin{picture}(75,10)
\put(10,6){\circle{3}}
\put(7,-4){\small 1}
\put(11.5,6){\line(1,0){10}}
\put(23,6){\circle{3}}
\put(20,-4){\small 2}
\put(24.5,6){\line(1,0){10}}
\put(36,6){\circle*{3}}
\put(33,-4){\small 3}
\put(49,6){\circle{3}}
\put(46,-4){\small 4}
\put(62,6){\circle*{3}}
\put(59,-4){\small 5}
\end{picture} Its adjacency matrix is a symmetric $5\times 5$-matrix $J$ with $J_{12}=1$, $J_{23}=1$.


\bigskip Before reformulating the definition of quasiplanar Wick products from~\cite{quasiI} in the present context, we need some more definitions.
We first extend the definition of regular chord diagrams~\cite{stoimenow}. 

{\defi \rm A 
dotted chord diagram is called {\em regular},
if for any two pairs of distinct points $(x,y)$ and $(w,z)$ in the diagram whose chords do not intersect, 
both $x$ and $y$ appear either on the right hand side or on the left hand side of both $z$ and $w$ on the arc.}

\medskip For example, the diagram \begin{picture}(50,10)
\put(10,0){\circle*{2}}
\put(20,0){\circle*{2}}
\put(30,0){\circle*{2}}
\put(40,0){\circle*{2}}
\qbezier(10,0)(25,16)(40,0)
\put(3,0){\vector(1,0){47}}
\end{picture} is
regular (and connected), while the diagram \begin{picture}(50,10)
\put(10,0){\circle*{2}}
\put(20,0){\circle*{2}}
\put(30,0){\circle*{2}}
\put(40,0){\circle*{2}}
\qbezier(10,0)(25,16)(40,0)
\qbezier(20,0)(25,11)(30,0)
\put(3,0){\vector(1,0){47}}
\end{picture} is not
regular.

{\defi \rm A dotted chord diagram is called {\em quasiplanar}, 
if for any two pairs of points $(x,y)$ and $(z,z)$ in the diagram, 
both $x$ and $y$ appear either on the right hand side or on the left hand side of $z$ on the arc. }

\medskip Observe that a diagram is quasiplanar if and only if its labelled intersection graph does not contain any edges between back and white vertices.
For example, the diagram 
\begin{picture}(95,10)
\put(10,0){\circle*{2}}
\put(20,0){\circle*{2}}
\put(30,0){\circle*{2}}
\put(40,0){\circle*{2}}
\put(50,0){\circle*{2}}
\put(60,0){\circle*{2}}
\put(70,0){\circle*{2}}
\put(80,0){\circle*{2}}
\qbezier(10,0)(20,16)(30,0)
\qbezier(20,0)(30,16)(40,0)
\qbezier(60,0)(65,11)(70,0)
\put(3,0){\vector(1,0){87}}
\end{picture} 
is quasiplanar (regular, not connected), while the diagram 
\begin{picture}(95,10)
\put(10,0){\circle*{2}}
\put(20,0){\circle*{2}}
\put(30,0){\circle*{2}}
\put(40,0){\circle*{2}}
\put(50,0){\circle*{2}}
\put(60,0){\circle*{2}}
\put(70,0){\circle*{2}}
\put(80,0){\circle*{2}}
\qbezier(10,0)(20,16)(30,0)
\qbezier(20,0)(35,16)(50,0)
\qbezier(60,0)(65,11)(70,0)
\put(3,0){\vector(1,0){87}}
\end{picture} 
is not quasiplanar (but still regular). 
Observe that any diagram without dots and any diagram without chords is quasiplanar. In particular, the diagram \begin{picture}(50,10)
\put(10,0){\circle*{2}}
\put(20,0){\circle*{2}}
\put(30,0){\circle*{2}}
\put(40,0){\circle*{2}}
\qbezier(10,0)(25,16)(40,0)
\qbezier(20,0)(25,11)(30,0)
\put(3,0){\vector(1,0){47}}
\end{picture} is quasiplanar (but not 
regular).

\medskip We will now consider regular quasiplanar diagrams, and denote by $V^{rq}_n$ the subspace of $V_n$ that is spanned by regular quasiplanar dotted chord diagrams of dot-degree $n$. We also use the notation $V^{rq}=\bigoplus V^{rq}_n$. 

{\rem \rm \label{regQuasi} It is not difficult to see that a product of connected quasiplanar diagrams is regular quasiplanar and that any regular quasiplanar diagram can be written in a unique way as a product of nontrivial connected quasiplanar diagrams. More generally, any regular diagram can be written in a unique way as a product of nontrivial connected  diagrams.}

\medskip For this reason, connected quasiplanar diagrams will turn out to be important.
Observe in particular, that a connected quasiplanar dotted chord diagram $D$ is either
the diagram \begin{picture}(20,10)
\put(5,3){\circle*{2}}
\put(0,3){\vector(1,0){15}}
\end{picture}
or it 
does not contain any dots. Unless it is of dot-degree 1, a connected quasiplanar dotted diagram therefore has even dot-degree. We will denote by $\mc D^{cq}_n$ the set of all connected quasiplanar dotted chord diagrams of dot-degree $n$, and by $\mc D^{cq}$ the set of all connected quasiplanar dotted chord diagrams.

\medskip In what follows, let us denote an (unlabelled) diagram without chords of dot-degree $n$ by $[n]$,
\[
[n]\ = \ 
\begin{picture}(75,10)
\put(10,0){\circle*{2}}
\put(20,0){\circle*{2}}
\put(35,0){\dots}
\put(60,0){\circle*{2}}
\put(70,0){\circle*{2}}
\put(3,0){\line(1,0){27}}
\put(55,0){\vector(1,0){27}}
\end{picture} 
\] 
with the convention $[0]=\emptyset$. 

{\defi \label{qplWick} \rm We define a linear map $\mc W : V \rightarrow V$, called the {\em quasiplanar Wick map}, by setting $\mc W(D)=0$ if $D$ contains a chord, and for chordless diagrams, we define $\mc W$   
inductively by $\mc W(\emptyset)= \emptyset$, and for $n\geq 1$, 
\[
\mc W([n])=\sum_{0 \neq m \leq n} (-1)^{m+1}\;f([m])\;
\mc W([n-m])
\qquad \mbox{ with } \qquad  f([m])=\sum_{D \in \mc D^{cq}_m} D 
\]
where the sum in the definition of $f$ runs over all connected quasiplanar diagrams of dot-degree $m$.  We call the image $\mc W([n])$ the $n$-fold {\em quasiplanar Wick product}. }

\medskip {\prop\label{WickContractions} \rm For any $n\geq 1$, we have 
\begin{equation}\label{eq:WickContractions}
\mc W([n])=\sum_{K=1}^n 
%
%
\sum_{D_i} (-1)^{n+K}\;D_1\cdots D_K
\end{equation}
where the sum runs over all nontrivial diagrams $D_i \in \mc D^{cq}$ whose dot-degrees are a partition of $n$, i.e. $\sum_{i=1}^K\deg D_i = n$.

}

\proof{The claim is almost obvious from the definition, although the formal proof turns out to appear complicated. Clearly, the claim is true for $n=1$, since $\mc W([1])=[1]$ which is the only connected quasiplanar diagram of dot-degree $1$.
Now, assume the claim to be true for $[n]$. Then by the induction hypothesis and inserting the definition for $f([2l])$, we have
\begin{eqnarray}
	\mc W([n+1]) &=&
		f([1]) \mc W([n]) - \sum_{l=1}^{\left[\frac{n+1}2\right]}
				f([2l])\,\mc W([n+1-2l])
\nonumber \\		&=& 
	[1] \; \sum_{K=1}^n \sum_{\tiny \begin{array}{c} D_i \in \mc D^{cq} \\ \sum\limits_{i=1}^K\deg D_i = n \end{array}} (-1)^{n+K}\;D_1\cdots D_K
\label{first}\\&&
 	- \sum_{l=1}^{\left[\frac{n+1}2\right]}
	\sum_{D_0 \in \mc D^{cq}_{2l}} \sum_{K=1}^{n+1-2l} 
	\hspace{-14ex} \sum_{\tiny \begin{array}{l} \hspace{25ex}
	D_i \in \mc D^{cq} 
	\\ \hspace{25ex}\sum\limits_{i=1}^K\deg D_i = {n+1-2l} \end{array}} \hspace{-12ex} (-1)^{n+1+K-2l}\;D_0 \, D_1\cdots D_K \quad 
\label{second}\end{eqnarray}
Now, line~(\ref{second}) in the above can be rewritten as 
\[
-\sum_{K=0}^{n-1} \sum_{D_i}
(-1)^{n+1+K}\;D_0 \, D_1\cdots D_K
\]
where the second sum runs over all nontrivial diagrams $D_0,D_1, \dots, D_K  \in \mc D^{cq}$ with $\sum_{i=0}^K\deg D_i = {n+1}$ and $\deg D_0 \geq 2$. Observe that for $n$ even, the sum actually starts with $K=1$, since for $K=0$, the second sum is empty (we would have $\deg D_0=n+1$, in contradiction with the fact that the dot-degree of $D_0$ has to be even). 
Shifting the summation index by one, we then find that line~(\ref{second}) is equal to  
\[
-\sum_{M=1}^{n} \sum_{D_i}
(-1)^{n+M}\; D_1\cdots D_M
\]
where the second sum now runs over all nontrivial diagrams $D_1, \dots, D_M  \in \mc D^{cq}$ with $\sum_{i=1}^M\deg D_i = {n+1}$ and $\deg D_1 \geq 2$. Observe that this sum can be extended to include the case $M=n+1$, since the sum over the diagrams is empty in that case (since the conditions $\deg D_1 \geq 2$ and $\deg D_1 + \dots +\deg D_{n+1}=n+1$ cannot be fulfilled simultanously).
Using that the diagram $[1]$ is in fact the sum over all diagrams $D_0 \in \mc D^{cq}$ with $\deg D_0=1$, and again shifting the summation index, we now rewrite line~(\ref{first}) in the above as follows
\[
\sum_{M=2}^{n+1} \sum_{D_i} (-1)^{n+M+1}\;D_1\cdots D_M
\]
where the second sum runs over all nontrivial diagrams $D_1, \dots, D_M \in \mc D^{cq}$ with $\sum_{i=1}^M\deg D_i = n +1$ and $\deg D_1=1$. This sum can in fact be extended to include $M=1$, since for $n\geq 1$,  the sum over the diagrams  is empty in this case anyway. Putting both sums together proves  the proposition.
}

\medskip 

\ex \label{W4} With the notation from above, we have
\begin{eqnarray*}
\mc W([4])  &=& [1] \; \mc W([3]) 
\ - \begin{picture}(35,10)
\put(10,0){\circle*{2}}
\put(20,0){\circle*{2}}
\qbezier(10,0)(15,16)(20,0)
\put(3,0){\vector(1,0){27}}
\end{picture} 
\; \mc W([2])
\ - 
\ 
\begin{picture}(55,10)
\put(10,0){\circle*{2}}
\put(20,0){\circle*{2}}
\put(30,0){\circle*{2}}
\put(40,0){\circle*{2}}
\qbezier(10,0)(20,16)(30,0)
\qbezier(20,0)(30,16)(40,0)
\put(3,0){\vector(1,0){47}}
\end{picture} 
\\ 
& = & 
\ \underbrace{\phantom{\int}\hspace{-4ex}
\begin{picture}(53,10)
\put(10,0){\circle*{2}}
\put(20,0){\circle*{2}}
\put(30,0){\circle*{2}}
\put(40,0){\circle*{2}}
\put(3,0){\vector(1,0){47}}
\end{picture} 
}_{\displaystyle K=4}
\; \underbrace{-\phantom{\int}\hspace{-2ex}
\begin{picture}(53,10)
\put(10,0){\circle*{2}}
\put(20,0){\circle*{2}}
\put(30,0){\circle*{2}}
\put(40,0){\circle*{2}}
\qbezier(10,0)(15,16)(20,0)
\put(3,0){\vector(1,0){47}}
\end{picture} 
-
\begin{picture}(55,10)
\put(10,0){\circle*{2}}
\put(20,0){\circle*{2}}
\put(30,0){\circle*{2}}
\put(40,0){\circle*{2}}
\qbezier(20,0)(25,16)(30,0)
\put(3,0){\vector(1,0){47}}
\end{picture} 
-
\begin{picture}(53,10)
\put(10,0){\circle*{2}}
\put(20,0){\circle*{2}}
\put(30,0){\circle*{2}}
\put(40,0){\circle*{2}}
\qbezier(30,0)(35,16)(40,0)
\put(3,0){\vector(1,0){47}}
\end{picture} 
}_{\displaystyle K=3} 
\\&&
 \underbrace{+\phantom{\int}\hspace{-2ex}
\begin{picture}(53,10)
\put(10,0){\circle*{2}}
\put(20,0){\circle*{2}}
\put(30,0){\circle*{2}}
\put(40,0){\circle*{2}}
\qbezier(10,0)(15,16)(20,0)
\qbezier(30,0)(35,16)(40,0)
\put(3,0){\vector(1,0){47}}
\end{picture} 
}_{\displaystyle K=2}
 \underbrace{-\phantom{\int}\hspace{-2ex}
\begin{picture}(53,10)
\put(10,0){\circle*{2}}
\put(20,0){\circle*{2}}
\put(30,0){\circle*{2}}
\put(40,0){\circle*{2}}
\qbezier(10,0)(20,16)(30,0)
\qbezier(20,0)(30,16)(40,0)
\put(3,0){\vector(1,0){47}}
\end{picture} 
}_{\displaystyle K=1}
\end{eqnarray*}

\medskip {\cor \label{EigenschaftenW} \rm The image of the quasiplanar Wick map $\mc W$ is contained in the space of regular quasiplanar diagrams $V^{rq}$. It is a projection, 
$\mc W \circ \mc W = \mc W$.
Moreover,  we have $\Delta\big(\,\mc W([n])\,\big)\;\Delta\big(\,\mc W([m])\,\big)=\Delta\big(\,\mc W([n])\, \mc W([m])\,\big)$.}

\proof{All claims are a consequence of the equality given by Proposition~\ref{WickContractions}.
The first claim follows immediately, since all terms on the right hand side of the equation are regular and quasiplanar. All but the last term in the sum on the right hand side (where $K=\deg([n])=n$) contain at least one chord, so the second claim follows from the definition of $\mc W$. The third claim is a consequence of the fact that $(V,\mu,\Delta)$ is a bialgebra.}

{\prop \label{qplWickthm0} \rm For the product of two quasiplanar Wick products $\mc W([n])$ and $\mc W([m])$, we find the relation
\begin{equation}\label{eq:qWickthm0}
\mc W([n]) \mc W([m]) = \mc W([n+m]) +  \sum_{K=1}^{n+m} \sum_{D_i} (-1)^{n+m+K+1}\;
D_1\cdots D_K
\end{equation}
where the sum runs over all nontrivial diagrams $D_1,\dots, D_K \in \mc D^{cq}$
such that $\sum_{i=1}^K\deg D_i = n+m$ but where no $1\leq S \leq n+m$ exists, such that $\sum_{i=1}^S\deg D_i = n$ and $\sum_{i=S+1}^K\deg D_i = m$.
}

\proof{We write down the expressions for $\mc W([n])$, $\mc W([m])$, and $\mc W([n+m])$ according to Proposition~\ref{WickContractions}. It is then not difficult to establish that all terms which appear in the product $\mc W([n])\mc W([m])$ also appear in $\mc W([n+m])$. The converse is not true; the diagrams that appear in $\mc W([n+m])$ but not in $\mc W([n])\mc W([m])$ are the regular diagrams which contain chords connecting some of the first $n$ points with some of the last $m-n+1$ points. These diagrams are of the form $D_1\cdots D_K$ with $D_i \in \mc D^{cq}$
such that $\sum_{i=1}^K\deg D_i = n+m$ but where no $1\leq S \leq n+m$ exists, such that $\sum_{i=1}^S\deg D_i = n$ and $\sum_{i=S+1}^K\deg D_i = m$. They each appear with  prefactor $(-1)^{n+m+K}$. Subtracting all such diagrams (with their prefactors) from $\mc W([n+m])$ therefore yields the product $\mc W([n])\mc W([m])$.
This proves the proposition.
}

\medskip 
The signs which appear 
in equation~(\ref{eq:qWickthm0}) above can also be given in terms of the number of connected diagrams with dot-degree strictly larger than 1. 
More generally, we have:

{\rem \label{signs} \rm For any product $D_1\cdots D_K$ with nontrivial diagrams $D_i \in \mc D^{cq}$, and $\sum \deg D_i = n$, we have 
\[
(-1)^{n+K}\;D_1\cdots D_K \ = \  (-1)^{d_2}\; D_1\cdots D_K
\]
where $d_2$ is the number of connected quasiplanar diagrams in $D_1\cdots D_K$ with dot-degree strictly greater than 1. 
To see that this is true, observe that $(-1)^{n+K}= (-1)^{n-K}$, and that $K=d_1+d_2$, where $d_1$ counts the number of connected diagrams of dot-degree 1. Now, any quasiplanar connected diagram of dot-degree $>1$ has even dot-degree, hence $n-d_1$ is even, and the claim follows.
}

\ex We have 
\begin{eqnarray*}
\mc W([2])  \; \mc W([2]) 
&=& \mc W([4]) \ 
\underbrace{+\phantom{\int}\hspace{-2ex}
\begin{picture}(53,10)
\put(10,0){\circle*{2}}
\put(20,0){\circle*{2}}
\put(30,0){\circle*{2}}
\put(40,0){\circle*{2}}
\qbezier(10,0)(20,16)(30,0)
\qbezier(20,0)(30,16)(40,0)
\put(3,0){\vector(1,0){47}}
\end{picture} 
}_{\displaystyle K=1}
\; \underbrace{+\phantom{\int}\hspace{-2ex}
\begin{picture}(53,10)
\put(10,0){\circle*{2}}
\put(20,0){\circle*{2}}
\put(30,0){\circle*{2}}
\put(40,0){\circle*{2}}
\qbezier(20,0)(25,16)(30,0)
\put(3,0){\vector(1,0){47}}
\end{picture} 
}_{\displaystyle K=3}
\end{eqnarray*}
and 
\begin{eqnarray*}
\mc W([2])  \; \mc W([3]) 
&=& \mc W([5]) \ 
 \underbrace{+\phantom{\int}\hspace{-2ex}
\begin{picture}(63,10)
\put(10,0){\circle*{2}}
\put(20,0){\circle*{2}}
\put(30,0){\circle*{2}}
\put(40,0){\circle*{2}}
\put(50,0){\circle*{2}}
\qbezier(20,0)(30,16)(40,0)
\qbezier(30,0)(40,16)(50,0)
\put(3,0){\vector(1,0){57}}
\end{picture} 
+
\begin{picture}(63,10)
\put(10,0){\circle*{2}}
\put(20,0){\circle*{2}}
\put(30,0){\circle*{2}}
\put(40,0){\circle*{2}}
\put(50,0){\circle*{2}}
\qbezier(10,0)(20,16)(30,0)
\qbezier(20,0)(30,16)(40,0)
\put(3,0){\vector(1,0){57}}
\end{picture} 
}_{\displaystyle K=2} \ 
\\&&
 \underbrace{-\phantom{\int}\hspace{-2ex}
\begin{picture}(63,10)
\put(10,0){\circle*{2}}
\put(20,0){\circle*{2}}
\put(30,0){\circle*{2}}
\put(40,0){\circle*{2}}
\put(50,0){\circle*{2}}
\qbezier(20,0)(25,16)(30,0)
\qbezier(40,0)(45,16)(50,0)
\put(3,0){\vector(1,0){57}}
\end{picture} 
}_{\displaystyle K=3} \ 
 \underbrace{+\phantom{\int}\hspace{-2ex}
\begin{picture}(63,10)
\put(10,0){\circle*{2}}
\put(20,0){\circle*{2}}
\put(30,0){\circle*{2}}
\put(40,0){\circle*{2}}
\put(50,0){\circle*{2}}
\qbezier(20,0)(25,16)(30,0)
\put(3,0){\vector(1,0){57}}
\end{picture} 
}_{\displaystyle K=4}
\end{eqnarray*}

\noindent 
For reasons inherent to quantum field theory, it is desirable to rewrite 
all terms on the right hand side of equation~(\ref{eq:qWickthm0}) in terms of the quasiplanar Wick map $\mc W$. This is achieved by first extending $\mc W$ to diagrams with chords. The image of this extension is contained in the vector space of quasiplanar diagrams, but in general no longer in the vector space of regular quasiplanar diagrams. The construction is not yet fully understood in algebraic terms and is not necessary to understand the results presented in the remaining sections of this paper.

By Proposition~\ref{WickContractions}, a dotted diagram without chords $[n]$ can obviously be rewritten as follows
\[
[n]=\mc W([n]) - \sum_{K=1}^{n-1} 
\sum_{D_i} (-1)^{n+K}\;D_1\cdots D_K
\]
where the sum runs over all nontrivial diagrams $D_i \in \mc D^{cq}$ with $\sum_{i=1}^K\deg D_i = n$. We would like to iterate this process of replacing dots by Wick products also in diagrams containing chords (such as the terms of the sum over $K$ in the above).

The map $\mc W$ itself cannot be used to that end, since it maps any diagram containing a chord to~0. We now define an extension $\mc W^\prime$ of $\mc W$ which is equal to $\mc W$ on diagrams without chords and on quasiplanar diagrams with chords acts as $\mc W$ on the diagram's dots while leaving the rest of the diagram unchanged. Observe that $\mc W^\prime$ takes values in the vector space of quasiplanar, but not necessarily regular diagrams, i.e. the image of $\mc W^\prime$ is not in general contained in $V^{rq}$.
We again call a term of the form $\mc W^\prime(D)$, where $D$ is a diagram, a quasiplanar Wick product. 

\ex We have 
\begin{eqnarray*}
\mc W^\prime(\begin{picture}(75,10)
\put(10,0){\circle*{2}}
\put(20,0){\circle*{2}}
\put(30,0){\circle*{2}}
\put(40,0){\circle*{2}}
\put(50,0){\circle*{2}}
\put(60,0){\circle*{2}}
\qbezier(40,0)(45,11)(50,0)
\put(3,0){\vector(1,0){67}}
\end{picture})
&=& \ \ 
\begin{picture}(75,10)
\put(10,0){\circle*{2}}
\put(20,0){\circle*{2}}
\put(30,0){\circle*{2}}
\put(40,0){\circle*{2}}
\put(50,0){\circle*{2}}
\put(60,0){\circle*{2}}
\qbezier(40,0)(45,11)(50,0)
\put(3,0){\vector(1,0){67}}
\end{picture}
\ - \ 
\begin{picture}(75,10)
\put(10,0){\circle*{2}}
\put(20,0){\circle*{2}}
\put(30,0){\circle*{2}}
\put(40,0){\circle*{2}}
\put(50,0){\circle*{2}}
\put(60,0){\circle*{2}}
\qbezier(10,0)(15,11)(20,0)
\qbezier(40,0)(45,11)(50,0)
\put(3,0){\vector(1,0){67}}
\end{picture}
\\
&&\!\! - \ 
\begin{picture}(75,10)
\put(10,0){\circle*{2}}
\put(20,0){\circle*{2}}
\put(30,0){\circle*{2}}
\put(40,0){\circle*{2}}
\put(50,0){\circle*{2}}
\put(60,0){\circle*{2}}
\qbezier(20,0)(25,11)(30,0)
\qbezier(40,0)(45,11)(50,0)
\put(3,0){\vector(1,0){67}}
\end{picture}
- \ 
\begin{picture}(75,10)
\put(10,0){\circle*{2}}
\put(20,0){\circle*{2}}
\put(30,0){\circle*{2}}
\put(40,0){\circle*{2}}
\put(50,0){\circle*{2}}
\put(60,0){\circle*{2}}
\qbezier(30,0)(45,20)(60,0)
\qbezier(40,0)(45,11)(50,0)
\put(3,0){\vector(1,0){67}}
\end{picture}
\\
&&\!\! + \ 
\begin{picture}(75,10)
\put(10,0){\circle*{2}}
\put(20,0){\circle*{2}}
\put(30,0){\circle*{2}}
\put(40,0){\circle*{2}}
\put(50,0){\circle*{2}}
\put(60,0){\circle*{2}}
\qbezier(10,0)(15,11)(20,0)
\qbezier(30,0)(45,20)(60,0)
\qbezier(40,0)(45,11)(50,0)
\put(3,0){\vector(1,0){67}}
\end{picture}
- \ 
\begin{picture}(75,10)
\put(10,0){\circle*{2}}
\put(20,0){\circle*{2}}
\put(30,0){\circle*{2}}
\put(40,0){\circle*{2}}
\put(50,0){\circle*{2}}
\put(60,0){\circle*{2}}
\qbezier(40,0)(45,11)(50,0)
\qbezier(10,0)(20,16)(30,0)
\qbezier(20,0)(39,19)(60,0)
\put(3,0){\vector(1,0){67}}
\end{picture}
\end{eqnarray*}
Compare this with the quasiplanar Wick product $\mc W([4])$ given on page~\pageref{W4}.

\medskip
Iterating the procedure of replacing products of dots by quasiplanar Wick products using $\mc W^\prime$, we can rewrite any quasiplanar dotted chord diagram $D$ in terms of elements of the image of $\mc W^\prime$ and of diagrams without dots. This means that the image of $\mc W^\prime$ together with diagrams without dots provides a basis for quasiplanar dotted chord diagrams. The  algebraic meaning of this, however,  remains to be understood.

\ex The diagram \begin{picture}(75,10)
\put(10,0){\circle*{2}}
\put(20,0){\circle*{2}}
\put(30,0){\circle*{2}}
\put(40,0){\circle*{2}}
\put(50,0){\circle*{2}}
\put(60,0){\circle*{2}}
\qbezier(40,0)(45,11)(50,0)
\put(3,0){\vector(1,0){67}}
\end{picture}
is equal to the sum \begin{eqnarray*}
&=& \ \ 
\mc W^\prime(\begin{picture}(75,20)
\put(10,0){\circle*{2}}
\put(20,0){\circle*{2}}
\put(30,0){\circle*{2}}
\put(40,0){\circle*{2}}
\put(50,0){\circle*{2}}
\put(60,0){\circle*{2}}
\qbezier(40,0)(45,11)(50,0)
\put(3,0){\vector(1,0){67}}
\end{picture})
\ - \ 
\begin{picture}(75,10)
\put(10,0){\circle*{2}}
\put(20,0){\circle*{2}}
\put(30,0){\circle*{2}}
\put(40,0){\circle*{2}}
\put(50,0){\circle*{2}}
\put(60,0){\circle*{2}}
\qbezier(10,0)(15,11)(20,0)
\qbezier(30,0)(45,20)(60,0)
\qbezier(40,0)(45,11)(50,0)
\put(3,0){\vector(1,0){67}}
\end{picture}
+ \ 
\begin{picture}(75,10)
\put(10,0){\circle*{2}}
\put(20,0){\circle*{2}}
\put(30,0){\circle*{2}}
\put(40,0){\circle*{2}}
\put(50,0){\circle*{2}}
\put(60,0){\circle*{2}}
\qbezier(40,0)(45,11)(50,0)
\qbezier(10,0)(20,16)(30,0)
\qbezier(20,0)(39,19)(60,0)
\put(3,0){\vector(1,0){67}}
\end{picture}
\\
&& \  + \ 
\begin{picture}(75,10)
\put(10,0){\circle*{2}}
\put(20,0){\circle*{2}}
\put(30,0){\circle*{2}}
\put(40,0){\circle*{2}}
\put(50,0){\circle*{2}}
\put(60,0){\circle*{2}}
\qbezier(10,0)(15,11)(20,0)
\qbezier(40,0)(45,11)(50,0)
\put(3,0){\vector(1,0){67}}
\end{picture}
\ + \ 
\begin{picture}(75,10)
\put(10,0){\circle*{2}}
\put(20,0){\circle*{2}}
\put(30,0){\circle*{2}}
\put(40,0){\circle*{2}}
\put(50,0){\circle*{2}}
\put(60,0){\circle*{2}}
\qbezier(20,0)(25,11)(30,0)
\qbezier(40,0)(45,11)(50,0)
\put(3,0){\vector(1,0){67}}
\end{picture}
\ + \ 
\begin{picture}(75,15)
\put(10,0){\circle*{2}}
\put(20,0){\circle*{2}}
\put(30,0){\circle*{2}}
\put(40,0){\circle*{2}}
\put(50,0){\circle*{2}}
\put(60,0){\circle*{2}}
\qbezier(30,0)(45,20)(60,0)
\qbezier(40,0)(45,11)(50,0)
\put(3,0){\vector(1,0){67}}
\end{picture}
\\
&=& \ \ 
\mc W^\prime(\begin{picture}(75,20)
\put(10,0){\circle*{2}}
\put(20,0){\circle*{2}}
\put(30,0){\circle*{2}}
\put(40,0){\circle*{2}}
\put(50,0){\circle*{2}}
\put(60,0){\circle*{2}}
\qbezier(40,0)(45,11)(50,0)
\put(3,0){\vector(1,0){67}}
\end{picture})
\ - \ 
\begin{picture}(75,10)
\put(10,0){\circle*{2}}
\put(20,0){\circle*{2}}
\put(30,0){\circle*{2}}
\put(40,0){\circle*{2}}
\put(50,0){\circle*{2}}
\put(60,0){\circle*{2}}
\qbezier(10,0)(15,11)(20,0)
\qbezier(30,0)(45,20)(60,0)
\qbezier(40,0)(45,11)(50,0)
\put(3,0){\vector(1,0){67}}
\end{picture}
+ \ 
\begin{picture}(75,10)
\put(10,0){\circle*{2}}
\put(20,0){\circle*{2}}
\put(30,0){\circle*{2}}
\put(40,0){\circle*{2}}
\put(50,0){\circle*{2}}
\put(60,0){\circle*{2}}
\qbezier(40,0)(45,11)(50,0)
\qbezier(10,0)(20,16)(30,0)
\qbezier(20,0)(39,19)(60,0)
\put(3,0){\vector(1,0){67}}
\end{picture}
\\
&& \!\!\! + \ 
\mc W^\prime(\begin{picture}(75,15)
\put(10,0){\circle*{2}}
\put(20,0){\circle*{2}}
\put(30,0){\circle*{2}}
\put(40,0){\circle*{2}}
\put(50,0){\circle*{2}}
\put(60,0){\circle*{2}}
\qbezier(10,0)(15,11)(20,0)
\qbezier(40,0)(45,11)(50,0)
\put(3,0){\vector(1,0){67}}
\end{picture})
\ - \ 
\begin{picture}(75,10)
\put(10,0){\circle*{2}}
\put(20,0){\circle*{2}}
\put(30,0){\circle*{2}}
\put(40,0){\circle*{2}}
\put(50,0){\circle*{2}}
\put(60,0){\circle*{2}}
\qbezier(10,0)(15,11)(20,0)
\qbezier(30,0)(45,20)(60,0)
\qbezier(40,0)(45,11)(50,0)
\put(3,0){\vector(1,0){67}}
\end{picture}
\\&& 
\!\!\! + \ 
\mc W^\prime(\begin{picture}(75,15)
\put(10,0){\circle*{2}}
\put(20,0){\circle*{2}}
\put(30,0){\circle*{2}}
\put(40,0){\circle*{2}}
\put(50,0){\circle*{2}}
\put(60,0){\circle*{2}}
\qbezier(20,0)(25,11)(30,0)
\qbezier(40,0)(45,11)(50,0)
\put(3,0){\vector(1,0){67}}
\end{picture})
\ - \ 
\begin{picture}(75,10)
\put(10,0){\circle*{2}}
\put(20,0){\circle*{2}}
\put(30,0){\circle*{2}}
\put(40,0){\circle*{2}}
\put(50,0){\circle*{2}}
\put(60,0){\circle*{2}}
\qbezier(20,0)(25,11)(30,0)
\qbezier(10,0)(35,22)(60,0)
\qbezier(40,0)(45,11)(50,0)
\put(3,0){\vector(1,0){67}}
\end{picture}
\\
&&\!\!\! + \ 
\mc W^\prime(\begin{picture}(75,15)
\put(10,0){\circle*{2}}
\put(20,0){\circle*{2}}
\put(30,0){\circle*{2}}
\put(40,0){\circle*{2}}
\put(50,0){\circle*{2}}
\put(60,0){\circle*{2}}
\qbezier(30,0)(45,20)(60,0)
\qbezier(40,0)(45,11)(50,0)
\put(3,0){\vector(1,0){67}}
\end{picture})
\ - \ 
\begin{picture}(75,15)
\put(10,0){\circle*{2}}
\put(20,0){\circle*{2}}
\put(30,0){\circle*{2}}
\put(40,0){\circle*{2}}
\put(50,0){\circle*{2}}
\put(60,0){\circle*{2}}
\qbezier(10,0)(15,11)(20,0)
\qbezier(30,0)(45,20)(60,0)
\qbezier(40,0)(45,11)(50,0)
\put(3,0){\vector(1,0){67}}
\end{picture}
\end{eqnarray*}

\medskip In particular, we can use the map $\mc W^\prime$ to rewrite all terms on the right hand side of equation~(\ref{eq:qWickthm0}) from Proposition~\ref{qplWickthm0} in terms of quasiplanar Wick products. The resulting equation is what is called the quasiplanar Wick theorem in~\cite{quasiI,bahnsDiss}.


\section{The shuffle Hopf algebra}\label{shuffleHopf}

\medskip By Corollary~\ref{EigenschaftenW}, the image of the quasiplanar Wick map $\mc W$ is contained in the space of quasiplanar regular diagrams $V^{rq}$. We will now consider a Hopf structure on $V^{rq}$ such that the quasiplanar Wick map can be understood as a convolution in this Hopf algebra. Incidentally, the ordinary Hopf algebra of chord diagrams does not seem to be natural in this context, but instead we have to consider the Hopf structure.

By remark~\ref{regQuasi}, any regular quasiplanar diagram can be uniquely written as a product of connected quasiplanar diagrams. We use this to endow  $V^{rq}$ with a commutative product, the so-called shuffle product $\mu_\#:V^{rq}\otimes V^{rq}\rightarrow V^{rq}$,
\[
\mu_\# (D \otimes D^\prime)= \sum_{\sigma \in Sh_{n,m}} D_{\sigma(1)} \cdots 
D_{\sigma(n+m)}
\]
where $D=D_1\cdots D_n$ and $D^\prime=D_{n+1}\cdots D_{n+m}$ with $D_i$ nontrivial, quasiplanar, and connected, and where $Sh_{n,m}=S_{n+m}/S_n\times S_m$ is the set of shuffle permutations, that is all elements $\sigma$ of $S_{n+m}$ which leave the order of the first $n$ elements and that of the last $m$ elements unchanged, i.e. 
$\sigma(1) < \sigma(2) <\dots \sigma(n)$ and $\sigma(n+1)<\sigma(n+2)<\dots < \sigma(n+m)$. 
We will usually write $a\# b$ for $\mu_\#(a\otimes b)$.
It is known that the shuffle product allows for a Hopf structure with the deconcatenation product as coproduct $\Delta_{dc}:V^{rq}\rightarrow V^{rq}\otimes V^{rq}$,
\[
\Delta_{dc}(D)=1\otimes D + D \otimes 1 + \sum_{k=1}^{n-1} D_1 \cdots D_k \otimes D_{k+1}\cdots D_{n}
\]
for $D=D_1\cdots D_n$  with $D_i$ nontrivial, quasiplanar and connected. Here, $1$ denotes the empty diagram $\emptyset=[0]$ which is the unit for the shuffle product. The counit is the map $\epsilon$ that is equal to 1 on the empty diagram and 0 elsewhere. The antipode of this Hopf algebra is $S(D_1 \cdots D_n)=(-1)^n D_n \cdots D_1$ where the diagrams $D_i$ are nontrivial, quasiplanar, and connected. Observe that $(V,\mu_\#,\Delta_{dc})$ is a commutative non-cocommutative graded connected Hopf algebra.

\ex We have
\begin{eqnarray*}
\begin{picture}(55,10)
\put(10,0){\circle*{2}}
\put(20,0){\circle*{2}}
\put(30,0){\circle*{2}}
\put(40,0){\circle*{2}}
\qbezier(10,0)(20,16)(30,0)
\qbezier(20,0)(30,16)(40,0)
\put(3,0){\vector(1,0){47}}
\end{picture} 
&\# &
\begin{picture}(55,10)
\put(10,0){\circle*{2}}
\put(20,0){\circle*{2}}
\put(30,0){\circle*{2}}
\put(40,0){\circle*{2}}
\qbezier(20,0)(25,11)(30,0)
\put(3,0){\vector(1,0){47}}
\end{picture} \ = 
\\ &&= \ \ 
\begin{picture}(95,10)
\put(10,0){\circle*{2}}
\put(20,0){\circle*{2}}
\put(30,0){\circle*{2}}
\put(40,0){\circle*{2}}
\put(50,0){\circle*{2}}
\put(60,0){\circle*{2}}
\put(70,0){\circle*{2}}
\put(80,0){\circle*{2}}
\qbezier(10,0)(20,16)(30,0)
\qbezier(20,0)(30,16)(40,0)
\qbezier(60,0)(65,11)(70,0)
\put(3,0){\vector(1,0){87}}
\end{picture}
\ + \ 
\begin{picture}(95,10)
\put(10,0){\circle*{2}}
\put(20,0){\circle*{2}}
\put(30,0){\circle*{2}}
\put(40,0){\circle*{2}}
\put(50,0){\circle*{2}}
\put(60,0){\circle*{2}}
\put(70,0){\circle*{2}}
\put(80,0){\circle*{2}}
\qbezier(20,0)(30,16)(40,0)
\qbezier(30,0)(40,16)(50,0)
\qbezier(60,0)(65,11)(70,0)
\put(3,0){\vector(1,0){87}}
\end{picture}
\\ 
&& \qquad  + \   
\begin{picture}(95,10)
\put(10,0){\circle*{2}}
\put(20,0){\circle*{2}}
\put(30,0){\circle*{2}}
\put(40,0){\circle*{2}}
\put(50,0){\circle*{2}}
\put(60,0){\circle*{2}}
\put(70,0){\circle*{2}}
\put(80,0){\circle*{2}}
\qbezier(20,0)(25,11)(30,0)
\qbezier(40,0)(50,16)(60,0)
\qbezier(50,0)(60,16)(70,0)
\put(3,0){\vector(1,0){87}}
\end{picture}
\ + \ 
\begin{picture}(95,10)
\put(10,0){\circle*{2}}
\put(20,0){\circle*{2}}
\put(30,0){\circle*{2}}
\put(40,0){\circle*{2}}
\put(50,0){\circle*{2}}
\put(60,0){\circle*{2}}
\put(70,0){\circle*{2}}
\put(80,0){\circle*{2}}
\qbezier(20,0)(25,11)(30,0)
\qbezier(50,0)(60,16)(70,0)
\qbezier(60,0)(70,16)(80,0)
\put(3,0){\vector(1,0){87}}
\end{picture}
\end{eqnarray*}
and 
\begin{eqnarray*}
\Delta_{dc}(\begin{picture}(95,10)
\put(10,0){\circle*{2}}
\put(20,0){\circle*{2}}
\put(30,0){\circle*{2}}
\put(40,0){\circle*{2}}
\put(50,0){\circle*{2}}
\put(60,0){\circle*{2}}
\put(70,0){\circle*{2}}
\put(80,0){\circle*{2}}
\qbezier(10,0)(20,16)(30,0)
\qbezier(20,0)(30,16)(40,0)
\qbezier(60,0)(65,11)(70,0)
\put(3,0){\vector(1,0){87}}
\end{picture}) \hspace{-26ex} && \hspace{26ex} \ = \\
&&  = \ 
1\otimes \begin{picture}(95,10)
\put(10,0){\circle*{2}}
\put(20,0){\circle*{2}}
\put(30,0){\circle*{2}}
\put(40,0){\circle*{2}}
\put(50,0){\circle*{2}}
\put(60,0){\circle*{2}}
\put(70,0){\circle*{2}}
\put(80,0){\circle*{2}}
\qbezier(10,0)(20,16)(30,0)
\qbezier(20,0)(30,16)(40,0)
\qbezier(60,0)(65,11)(70,0)
\put(3,0){\vector(1,0){87}}
\end{picture}
 + 
 \begin{picture}(95,10)
\put(10,0){\circle*{2}}
\put(20,0){\circle*{2}}
\put(30,0){\circle*{2}}
\put(40,0){\circle*{2}}
\put(50,0){\circle*{2}}
\put(60,0){\circle*{2}}
\put(70,0){\circle*{2}}
\put(80,0){\circle*{2}}
\qbezier(10,0)(20,16)(30,0)
\qbezier(20,0)(30,16)(40,0)
\qbezier(60,0)(65,11)(70,0)
\put(3,0){\vector(1,0){87}}
\end{picture} \otimes 1 
\\
&&\ \ \  + \ 
\begin{picture}(55,10)
\put(10,0){\circle*{2}}
\put(20,0){\circle*{2}}
\put(30,0){\circle*{2}}
\put(40,0){\circle*{2}}
\qbezier(10,0)(20,16)(30,0)
\qbezier(20,0)(30,16)(40,0)
\put(3,0){\vector(1,0){47}}
\end{picture} \otimes  
\begin{picture}(55,10)
\put(10,0){\circle*{2}}
\put(20,0){\circle*{2}}
\put(30,0){\circle*{2}}
\put(40,0){\circle*{2}}
\qbezier(20,0)(25,11)(30,0)
\put(3,0){\vector(1,0){47}}
\end{picture}
\ + \ 
\begin{picture}(65,10)
\put(10,0){\circle*{2}}
\put(20,0){\circle*{2}}
\put(30,0){\circle*{2}}
\put(40,0){\circle*{2}}
\put(50,0){\circle*{2}}
\qbezier(10,0)(20,16)(30,0)
\qbezier(20,0)(30,16)(40,0)
\put(3,0){\vector(1,0){57}}
\end{picture} \otimes  
\begin{picture}(34,10)
\put(10,0){\circle*{2}}
\put(20,0){\circle*{2}}
\put(30,0){\circle*{2}}
\qbezier(10,0)(15,11)(20,0)
\put(3,0){\vector(1,0){37}}
\end{picture}
\\
&&\ \ \  + \ 
\begin{picture}(85,10)
\put(10,0){\circle*{2}}
\put(20,0){\circle*{2}}
\put(30,0){\circle*{2}}
\put(40,0){\circle*{2}}
\put(50,0){\circle*{2}}
\put(60,0){\circle*{2}}
\put(70,0){\circle*{2}}
\qbezier(10,0)(20,16)(30,0)
\qbezier(20,0)(30,16)(40,0)
\qbezier(60,0)(65,11)(70,0)
\put(3,0){\vector(1,0){77}}
\end{picture} \otimes  
\begin{picture}(25,10)
\put(10,0){\circle*{2}}
\put(3,0){\vector(1,0){17}}
\end{picture}
\end{eqnarray*}

\medskip By definition, the primitive elements of this Hopf algebra are the connected quasiplanar diagrams.


{\prop \rm Let $h:V\rightarrow V$ be the map defined by $h([0])=1$,  $h(D)=0$ for any diagram $D$ containing chords and for any diagram of the form $D=[2k+1]$, and let
\[
h([n])=\sum_{K=1}^{n}(-1)^K \sum_{D_i} D_1\cdots D_K
\]
for even $n$, where the second sum runs over all connected quasiplanar diagrams of dot-degree strictly greater than 1 with $\sum \deg D_i = n$. Then the quasiplanar Wick product on diagrams without chords is the convolution of the identity map with $h$ with respect to the shuffle Hopf algebra,
\[
\mc W([n]) = \id \star h \;([n])  = \mu_\# \circ (\id \otimes h) \circ \Delta_{dc} \;([n]) 
\]
}

\proof{By the definition of $\Delta_{dc}$, and observing that $h([0])=1$, we find 
\begin{eqnarray*}
\mu_\# \circ (\id \otimes h) \circ \Delta_{dc} \;([n]) 
&=& \mu_\# \Big(\sum_{r=0}^n [r] \otimes h([n-r])\,\Big)
\\
&=& \mu_\# \Big(1\otimes h([n]) + [n]\otimes 1 +  \sum_{r=1}^{n-1} [r]\otimes h([n-r]) 
\,\Big)
\end{eqnarray*}
Now, inserting the definition of $h$, we find for the third term 
\begin{eqnarray*}
\sum_{r=1}^{n-1} [r]\otimes h([n-r]) &=& 
\sum_{L=1}^{n-1} \sum_{D_i}  (-1)^{d_2}\; D_1 \dots D_r \otimes D_{r+1}\cdots D_L
\end{eqnarray*}
Here, the sum runs over all nontrivial quasiplanar connected diagrams $D_1,\dots, D_L$ with $\sum \deg D_i =n$ where at least one and at most $n-1$ diagrams are of dot-degree~1,  $d_2$ denotes the number of  
diagrams  with dot-degree $\geq 2$, and moreover, all diagrams $D_i$ of dot-degree $1$ are in the tensor product's  first entry. Observe that this last  condition fixes the value of $r$ in the above. 

Application of the shuffle product $\mu_\#$ then yields $h([n])+[n]$ for the first two terms and for the sum above it distributes the diagrams of dot-degree 1, i.e. the dots $[1]$ from the tensor product's first entry, in all possible orders between the connected components of higher dot-degree. 
This yields all the terms that appear on the right hand side of equation~(\ref{eq:WickContractions}) in Proposition~\ref{WickContractions}, and by remark~\ref{signs}, also the signs $(-1)^{d_2}$ match those appearing in equation~(\ref{eq:WickContractions}).  This proves the proposition.}

\medskip \ex For the diagram $[4]$, we have indeed  
\begin{eqnarray*}
\mu_\#(\id \otimes h)\Delta_{dc}\;([4]) &=&
\mu_\#(\id \otimes h)(1\otimes 
\begin{picture}(55,10)
\put(10,0){\circle*{2}}
\put(20,0){\circle*{2}}
\put(30,0){\circle*{2}}
\put(40,0){\circle*{2}}
\put(3,0){\vector(1,0){47}}
\end{picture} 
+
\begin{picture}(25,10)
\put(10,0){\circle*{2}}
\put(3,0){\vector(1,0){17}}
\end{picture} 
\otimes 
\begin{picture}(45,10)
\put(10,0){\circle*{2}}
\put(20,0){\circle*{2}}
\put(30,0){\circle*{2}}
\put(3,0){\vector(1,0){37}}
\end{picture} 
+
\\
&&\phantom{\mu_\#(\id \otimes h)(}
+
\begin{picture}(35,10)
\put(10,0){\circle*{2}}
\put(20,0){\circle*{2}}
\put(3,0){\vector(1,0){27}}
\end{picture} 
\otimes
\begin{picture}(35,10)
\put(10,0){\circle*{2}}
\put(20,0){\circle*{2}}
\put(3,0){\vector(1,0){27}}
\end{picture} 
+
\\
&&\phantom{\mu_\#(\id \otimes h)(}
+ \begin{picture}(45,10)
\put(10,0){\circle*{2}}
\put(20,0){\circle*{2}}
\put(30,0){\circle*{2}}
\put(3,0){\vector(1,0){37}}
\end{picture} 
\otimes 
\begin{picture}(25,10)
\put(10,0){\circle*{2}}
\put(3,0){\vector(1,0){17}}
\end{picture} 
+
\begin{picture}(55,10)
\put(10,0){\circle*{2}}
\put(20,0){\circle*{2}}
\put(30,0){\circle*{2}}
\put(40,0){\circle*{2}}
\put(3,0){\vector(1,0){47}}
\end{picture} 
\otimes 1)
\\
&=&\mu_\#(
-\ 1\otimes 
\begin{picture}(55,15)
\put(10,0){\circle*{2}}
\put(20,0){\circle*{2}}
\put(30,0){\circle*{2}}
\put(40,0){\circle*{2}}
\qbezier(10,0)(20,16)(30,0)
\qbezier(20,0)(30,16)(40,0)
\put(3,0){\vector(1,0){47}}
\end{picture} 
+ 1\otimes 
\begin{picture}(55,10)
\put(10,0){\circle*{2}}
\put(20,0){\circle*{2}}
\put(30,0){\circle*{2}}
\put(40,0){\circle*{2}}
\qbezier(10,0)(15,11)(20,0)
\qbezier(30,0)(35,11)(40,0)
\put(3,0){\vector(1,0){47}}
\end{picture} 
+\ 
0
\\
&&\phantom{\mu_\#(}
- 
\begin{picture}(35,10)
\put(10,0){\circle*{2}}
\put(20,0){\circle*{2}}
\put(3,0){\vector(1,0){27}}
\end{picture} 
\otimes
\begin{picture}(35,10)
\put(10,0){\circle*{2}}
\put(20,0){\circle*{2}}
\qbezier(10,0)(15,11)(20,0)
\put(3,0){\vector(1,0){27}}
\end{picture} 
+
\
+ \ 
0 
\ + 
\ 
\begin{picture}(55,10)
\put(10,0){\circle*{2}}
\put(20,0){\circle*{2}}
\put(30,0){\circle*{2}}
\put(40,0){\circle*{2}}
\put(3,0){\vector(1,0){47}}
\end{picture} 
\otimes 1)
\\
&=& \mc W([4]) \begin{picture}(55,15)
\end{picture} 
\end{eqnarray*}

\medskip \noindent There is reason to hope that the extension $\mc W^\prime$ of $\mc W$ can be understood in algebraic terms in a similar manner, and that an algebraic version of the quasiplanar Wick theorem can be given.


\medskip

\section{Relation to quantum fields on the noncommutative Minkowski space}\label{qft}

\medskip 

Let me now recall from~\cite{quasiI,bahnsDiss} how the graphical language of the previous sections encodes certain operator valued distributions of quantum field theory on the noncommutative Minkowski space.
The {\em noncommutative Minkowski space} 
\[
M_\theta := \C \langle q_0,q_1,q_2,q_3\rangle /R_\theta
\]
is a quotient of the free algebra of 4 generators $\C \langle q_0,q_1,q_2,q_3\rangle$,
where $R_\theta$ is the ideal defined by $q_\mu q_\nu-q_\nu q_\mu - i \theta_{\mu\nu}\,I$ for $\mu,\nu\in\{0,1,2,3\}$ with a {\em nondegenerate} antisymmetric matrix $(\theta_{\mu\nu}) \in M(4\times 4,\R)$. 
In order to define quantum fields on $M_\theta$, it is convenient to work with the corresponding Weyl algebra generated by the Weyl operators
\[
e^{ikq}\ ,  \quad \mbox{ where } kq= \ts{\sum\limits_{\mu=0}^3} k^\mu \, q_\mu\ \mbox{ with }  k\in \R^4\ , \  k^\mu=\ts{\sum\limits_{\nu=0}^3} k_\nu \,\eta^{\nu\mu} \ ,  \ \eta={\rm diag}(+,-,-,-)
\]
such that for $p,k\in \R^4$,
\[
e^{ikq}\,e^{ipq} = e^{-\frac i 2 k\theta p}\, e^{i(p+k)q}\qquad 
\mbox{ with } k\theta p= \ts{\sum\limits_{\mu,\nu =0}^3} k^\mu \,\theta_{\mu\nu}\, p^\nu 
\]
The signs in the symmetric form $\eta$ above (the Minkowski metric) are at this point merely a convention with no important consequences. The signature of $\eta$ will, however, play a decisive role when the partial differential operators that are relevant in field theory are considered. In fact, questions of renormalization substantially depend on the signature of $\eta$, as I will show elsewhere~\cite{bahnsIR}.


Consider the operator valued distribution  $\varphi$ given by the free massive scalar 
real Klein Gordon field. Let $\omega$ be a state (i.e. a positive linear functional) on the Weyl algebra, and let $\psi_{\omega}$ denote its associated
Wigner function whose Fourier transform is 
$\hat{\psi}_{\omega}(k)=\omega(e^{ikq})$. Then the free massive scalar real Klein Gordon field $\phi$ on quantum spacetime $\mathcal E$ is defined as
an affine functional on a dense set of the state space of the Weyl algebra,  with values in the endomorphisms of a dense subset of Fock space, see~\cite{quasiII},
by the equation
 \begin{equation}
  \label{eq:phi}
  \phi(\omega)=\varphi(\psi_{\omega}) \ .
\end{equation}
Let $n>0$, let $f$ be a Schwartz function on $\R^{4n}$, let   $\widehat{\psi_{\omega}^n}(k_1,\dots, k_n):=\omega(\prod_{j=1}^n e^{ik_j q})$, and let   
$\times$ denote the convolution. Then the {\em regularized power} of $\phi$  is defined by
\begin{equation}
  \label{eq:product2}
  \phi_f^n(\omega)=\varphi^{\otimes n}(\psi_{\omega}^n\times f) \ .
\end{equation}
where
$\varphi^{\otimes n}$ is the operator  valued distribution in $n$ variables,
formally defined by its integral kernel
$ \varphi^{\otimes n}(x_1,\dots, x_n)=\prod_{i=1}^n\varphi(x_i) 
$ as usual. 
In our graphical language, the regularized power $\phi_f^n$ is symbolized by a dotted diagram of dot-degree $n$ without chords,
\[
\phi_f^n \qquad \leftrightarrow \qquad [n]\ = \ 
\begin{picture}(75,10)
\put(10,0){\circle*{2}}
\put(20,0){\circle*{2}}
\put(35,0){\dots}
\put(60,0){\circle*{2}}
\put(70,0){\circle*{2}}
\put(3,0){\line(1,0){27}}
\put(55,0){\vector(1,0){27}}
\end{picture} 
\] 
We have defined a renormalization procedure~\cite{quasiI,bahnsDiss} for such powers of fields, based on a certain notion of locality, which has lead to the definition of quasiplanar Wick products. 
In~\cite{quasiII} we will complete the proof that this procedure leads to operator valued distributions that are still well-defined when the Schwartz function $f$ is replaced by a $\delta$-distribution. The functional analytic details however, are not our concern here.

Instead, I will take the definition of quasiplanar Wick products for granted and merely recall how the distributions which appear in that definition can be symbolized in terms of dotted chord diagrams. 

Let $D$ be a dotted chord diagram whose labelled intersection graph $G$ has adjacency matrix $J$. Let $N$ denote the set of labels in $G$, let $U\subset N$ and $A\subset N$ denote the subsets which label the black vertices (dots) and the white vertices (chords) in $G$, respectively. Let $f$ be a Schwartz function $f \in \mc S(\R^{4n})$ which is {\em symmetric} in its arguments, $f(x_1,\dots, x_n)=f(x_{\pi(1)},\dots, x_{\pi(n)})$ for any $\pi \in S_n$. Then the labelled intersection graph $G$ of $D$ determines a Schwartz function  $f_G \in \mc S(\R^{4u})$ where $u=|U|$ by setting its Fourier transform to
\[
  \widehat{f_G}(k_U)=\int \prod_{a\in A} d\mu(k_a)\ 
  \exp\big(-i\sum_{s<t \ \in N}  J_{st} \, k_s\theta k_t\; \big)
  \ 
 \hat{f}(k_A,-k_A,k_U) 
\]
with the Lorentz invariant measure on the positive mass shell $d\mu(k)$, and where for an index set $I$, the symbol $k_I$ abbreviates the tuple $(k_i)_{i \in I}$ and where each $k_i$ is an element of $\R^4$. 
The explicit form of these integrals is $\int d\mu(k) f(k) = \int d^3p\; \sqrt{p^2+m^2}^{\,-1}\;f(\sqrt{p^2+m^2}\,, p)$ for $f\in \mc S(\R^4)$, $p\in\R^3$. We extend the correspondence $D \leftrightarrow f_G$ by linearity.

The operator valued distribution corresponding to $G$, evaluated in a testfunction $f$,  is then given as follows
\[
\phi^{u}_{f_G}(\omega)  \qquad \mbox{ which  by  (\ref{eq:product2}) is equal to }
 \   \varphi^{\otimes u}(\psi_{\omega}^{u}\times f_G) \ .
\]
Observe that the socalled twisting $  \exp\big(-i\sum_{s<t \ \in N}  J_{st} \, k_s\theta k_t\; \big)$ contains products $k_u \theta k_a$ or $k_a \theta k_u$ where $u\in U$ and $a \in A$, if and only if the dotted diagram is not quasiplanar.

Translating all graphs from Definition~\ref{qplWick} according to the above indeed yields the definition of quasiplanar Wick products from~\cite{quasiI,bahnsDiss}. The quasiplanarity of all terms on the right hand side of the equation in Proposition~\ref{WickContractions} means that the corresponding distributions fulfill the locality condition requested in~\cite{quasiI,bahnsDiss}. 

Note that in the investigation of quantum fields on the noncommutative Minkowski space, we also considered non-regular diagrams. It is known that once we admit for non-regular diagrams, a chord diagram is no longer uniquely determined by its labelled intersection graph. For example,  the labelled intersection graph of 
\begin{picture}(50,10)
\put(10,0){\circle*{2}}
\put(20,0){\circle*{2}}
\put(30,0){\circle*{2}}
\put(40,0){\circle*{2}}
\qbezier(10,0)(25,16)(40,0)
\qbezier(20,0)(25,11)(30,0)
\put(3,0){\vector(1,0){47}}
\end{picture} is the same as that of 
\begin{picture}(50,10)
\put(10,0){\circle*{2}}
\put(20,0){\circle*{2}}
\put(30,0){\circle*{2}}
\put(40,0){\circle*{2}}
\qbezier(10,0)(15,11)(20,0)
\qbezier(30,0)(35,11)(40,0)
\put(3,0){\vector(1,0){47}}
\end{picture} To distinguish such diagrams, we used a slightly more complicated definition of $f_G$ than the one given above. It is based on the use of arbitrary Schwartz functions $f$ (which are in general not symmetric under reordering their arguments), 
where for $a \in A$, the positions of the arguments $\pm k_a$ in $\hat f$ encode the starting and end point of the chords. For instance, for the two examples above we would have 
\[
\int d\mu(k_1)d\mu(k_2)\ 
 \hat f(k_1,k_2,-k_2,-k_1)
\quad
\mbox{ and } \int d\mu(k_1)d\mu(k_2)\ 
\hat f(k_1,-k_1,k_2,k_2) \ ,
\]
respectively. In this manner, the diagram itself, 
not only its labelled intersection graph is encoded, 
\[
D \leftrightarrow f_D
\]
and in particular, non-regular graphs can be distinguished. This more general definition coincides with the one given above when we consider Schwartz functions that are symmetric under rearranging their arguments.


\section{Knots}

\medskip By a theorem of Kontsevich, universal properties of knot invariants of finite degree are captured in terms of the vector space of chord diagrams modulo certain ideals given by the so-called framing independence and the 4T relation. An important tool in this framework are weight systems.

A weight system of degree $m$ is a  linear functional on the quotient of chord diagrams of degree $m$ modulo the framing independence and the 4T relation.
It is natural to ask whether in the special special case of diagrams without dots, the correspondence $D \leftrightarrow f_D$, extended by linearity, defines a weight system.

Unfortunately, this does not seem to be the case, as the correspondence is not well-defined on the quotient, as the 4T relation apparently cannot be satisfied, even for modifications of the correspondence above. Let us recall that the framing independence relation is that an arbitrary chord diagram containing an isolated chord (i.e. one that does not intersect any other chord) is  0, and that the 4T relation in terms of diagrams reads
\begin{eqnarray*}
\begin{picture}(65,20)
\put(12,0){\circle*{2}}
\put(17,0){\circle*{2}}
\put(10,-4){\line(1,0){7}}
\put(32,0){\circle*{2}}
\put(47,0){\circle*{2}}
\qbezier(12,0)(30,30)(47,0)
\qbezier(17,0)(24.5,16)(32,0)
\put(3,0){\vector(1,0){59}}
\end{picture} 
-
\begin{picture}(65,20)
\put(12,0){\circle*{2}}
\put(17,0){\circle*{2}}
\put(10,-4){\line(1,0){7}}
\put(32,0){\circle*{2}}
\put(47,0){\circle*{2}}
\qbezier(17,0)(32,30)(47,0)
\qbezier(12,0)(22,16)(32,0)
\put(3,0){\vector(1,0){59}}
\end{picture} 
&=&
\begin{picture}(65,20)
\put(12,0){\circle*{2}}
\put(27,0){\circle*{2}}
\put(32,0){\circle*{2}}
\put(26,-4){\line(1,0){7}}
\put(47,0){\circle*{2}}
\qbezier(12,0)(19.5,16)(27,0)
\qbezier(32,0)(39.5,16)(47,0)
\put(3,0){\vector(1,0){59}}
\end{picture} 
-
\begin{picture}(65,20)
\put(12,0){\circle*{2}}
\put(27,0){\circle*{2}}
\put(32,0){\circle*{2}}
\put(26,-4){\line(1,0){7}}
\put(47,0){\circle*{2}}
\qbezier(12,0)(22,16)(32,0)
\qbezier(27,0)(37,16)(47,0)
\put(3,0){\vector(1,0){59}}
\end{picture} 
\\&=&
\begin{picture}(65,20)
\put(12,0){\circle*{2}}
\put(27,0){\circle*{2}}
\put(42,0){\circle*{2}}
\put(47,0){\circle*{2}}
\put(40,-4){\line(1,0){7}}
\qbezier(12,0)(30,30)(47,0)
\qbezier(27,0)(35.5,16)(42,0)
\put(3,0){\vector(1,0){59}}
\end{picture} 
-
\begin{picture}(65,20)
\put(12,0){\circle*{2}}
\put(27,0){\circle*{2}}
\put(42,0){\circle*{2}}
\put(47,0){\circle*{2}}
\put(40,-4){\line(1,0){7}}
\qbezier(12,0)(27,30)(42,0)
\qbezier(27,0)(37.5,16)(47,0)
\put(3,0){\vector(1,0){59}}
\end{picture} 
\end{eqnarray*}
where arbitrary chords may be added in all 6 diagrams in the same positions in such a way that the two points which are underlined remain neighbours.

We now consider the framing independence relation.
Let $D$ be a diagram (without dots) of degree $n$ (i.e. of dot-degree $2n$), let $J$ denote the adjacency matrix of the corresponding labelled intersection graph $G$. Let $f \in \mc S(\R^{3n}) $ be symmetric under reordering its arguments, and consider the map 
\begin{equation}\label{WeightCorresp}
G \quad \rightarrow  \quad 
\int \prod_{i=1}^n d\mu(p_i)\  
  \exp\big(-i\sum_{s<t}  J_{st} \, p_s\theta p_t\; \big)
  \ 
f({\bf p}_1, \dots, {\bf p}_n) 
\end{equation}
where $p=(\sqrt{{\bf p}^2+m^2}\,, {\bf p})$.
The right hand side is very similar to the definition of $f_G$, the sole difference being that we use a function $f({\bf k})$ in place of $\hat f(k,-k)$ where $k=(\sqrt{{\bf k}^2+m^2}\,, {\bf k})$. Now, let $j$ be the label of an isolated chord in $D$. By definition, we have $J_{jk}=0$ and $J_{kj}=0$ for all $k$. Hence, the integral above is zero on diagrams with isolated chords as desired provided we choose a Schwartz function $f$ with total integral 0.

We now try to implement the 4T relation as well. To that end, we use 
Lemma XX.3.1. from~\cite{kassel}:
Let $D$ be a diagram (without dots) containing at least 2 chords, let $x$ be the point on the arc that is furthest to the left, let $(x,y)$ be the corresponding chord. Consider a point $x^\prime$ to the right hand side of $D$ on the arc. Then the diagram $D^\prime=(D\setminus (x,y)) \cup (y,x^\prime)$ and $D$ coincide modulo the 4T relation.

In order to define a weight system, the correspondence (\ref{WeightCorresp}) must therefore yield the same result for any two diagrams 
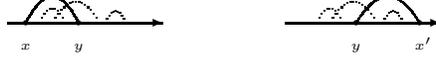
\[
\begin{picture}(65,30)
\put(10,10){\circle*{2}}
\put(8,0){{\tiny $x$}}
\put(30,10){\circle*{2}}
\put(28,0){{\tiny $y$}}
\qbezier(10,10)(20,30)(30,10)
\qbezier[16](20,10)(29,26)(38,10)
\qbezier[10](15,10)(20,21)(25,10)
\qbezier[10](40,10)(44,19)(48,10)
\put(3,10){\vector(1,0){59}}
\end{picture} 
 \qquad \qquad 
\begin{picture}(65,30)
\put(30,10){\circle*{2}}
\put(28,0){{\tiny $y$}}
\put(54,10){\circle*{2}}
\put(52,0){{\tiny $x^\prime$}}
\qbezier[16](20,10)(29,26)(38,10)
\qbezier[10](15,10)(20,21)(25,10)
\qbezier[10](40,10)(44,19)(48,10)
\qbezier(30,10)(42,30)(54,10)
\put(3,10){\vector(1,0){59}}
\end{picture} 
\]
with an arbitrary number of additional chords placed in the same positions in the two diagrams, as indicated by the dotted chords.

Now label the first of these diagrams starting with the label $0$. Then the map~(\ref{WeightCorresp}) gives the integral
\begin{eqnarray*}
&&\int d\mu(p_0) \prod_{i=1}^{n-1} d\mu(p_i)\  
  \exp\big(-i\sum_{0<t\leq n-1}  J_{0t} \, p_0\theta p_t\; \big)\, \cdot 
\\ && \qquad \qquad \cdot  \ 
  \exp\big(-i\sum_{0<s<t\leq n-1}  J_{st} \, p_s\theta p_t\; \big)
 \ f({\bf p}_0,\dots, {\bf p}_{n-1}) 
\end{eqnarray*}
Labelling the second diagram with labels $1,\dots,n$, on the other hand yields the integral
\[
\int \prod_{i=1}^{n} d\mu(p_i)\  
  \exp\big(-i\hspace{-2ex}\sum_{0<s<t\leq n-1} \hspace{-2ex} J_{st} \, p_s\theta p_t\; \big)
  \exp\big(-i\sum_{0<t\leq n-1}  J_{tn} \, p_t\theta p_n\; \big)\,
  \ 
 f({\bf p}_1,\dots, {\bf p}_n) 
\]
As the matrix $\theta$ is antisymmetric, one might think of changing the variables ${\bf p}_n=-{\bf p}_0$, but this would leave a wrong sign in the exponentials which involve the component $\sqrt{{\bf p}_n^2+m^2}$ of $p_n=(\sqrt{{\bf p}_n^2+m^2}\,, {\bf p}_n)$. We conclude that the map~(\ref{WeightCorresp}) cannot well defined on the quotient given by the 4T relation.

One might therefore try to slightly change the map (\ref{WeightCorresp}) by considering integrals over $\R^{2n}$ or $\R^{4n}$, instead of integrals over the mass shell\footnote{Observe that an antisymmetric matrix has even rank, so even if $\theta$ has maximal rank, a twisting $\exp(-i k\theta p)$ would be independent of one component of $k\in \R^{2k+1}$.},
\begin{equation}\label{WeightCorresp2D}
G \quad \rightarrow  \quad 
\int \prod_{i=1}^n d^4 p_i\  
  \exp\big(-i\sum_{s<t}  J_{st} \, p_s\theta p_t\; \big)
  \ 
f(p_1, \dots, p_n) 
\end{equation}
where $f\in \mc S(\R^{4n})$. In order to accomodate the framing independence relation, we would again ask that the total integrals $\int d p_i f(p_1,\dots, p_n)$ vanish. 

The attempt to also implement the 4T relation however, leaves no other possibility than to choose $f=0$, thus producing the trivial weight system. 
To see this, observe that according to the discussion of the map (\ref{WeightCorresp}) this would require that 
\[
f(p_1,\dots,p_{n-1},-p_0) \ = \ f(p_0,p_1,\dots,p_{n-1}) 
\]
in contradiction with the framing independence relation.

A last loophole might be to consider the map that is 0 on diagrams containing isolated chords and is otherwise given by (\ref{WeightCorresp2D}). However, even in this case, the 4T relation does not hold. To see this, suffice to calculate the twisting for 
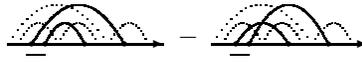
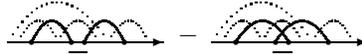
\[
\begin{picture}(65,20)
\put(12,0){\circle*{2}}
\put(17,0){\circle*{2}}
\put(10,-4){\line(1,0){7}}
\put(32,0){\circle*{2}}
\put(47,0){\circle*{2}}
\qbezier[15](7,0)(14.5,16)(22,0)
\qbezier[20](3,0)(21,30)(41,0)
\qbezier[15](22,0)(29.5,16)(37,0)
\qbezier[15](41,0)(49.5,16)(56,0)
\qbezier(12,0)(30,30)(47,0)
\qbezier(17,0)(24.5,16)(32,0)
\put(3,0){\vector(1,0){59}}
\end{picture} 
-
\begin{picture}(65,20)
\put(12,0){\circle*{2}}
\put(17,0){\circle*{2}}
\put(10,-4){\line(1,0){7}}
\put(32,0){\circle*{2}}
\put(47,0){\circle*{2}}
\qbezier[15](7,0)(14.5,16)(22,0)
\qbezier[20](3,0)(21,30)(41,0)
\qbezier[15](22,0)(29.5,16)(37,0)
\qbezier[15](41,0)(49.5,16)(56,0)
\qbezier(17,0)(32,30)(47,0)
\qbezier(12,0)(22,16)(32,0)
\put(3,0){\vector(1,0){59}}
\end{picture} 
\]
and
\[
\begin{picture}(65,20)
\put(12,0){\circle*{2}}
\put(27,0){\circle*{2}}
\put(32,0){\circle*{2}}
\put(26,-4){\line(1,0){7}}
\put(47,0){\circle*{2}}
\qbezier[15](7,0)(14.5,16)(22,0)
\qbezier[20](3,0)(21,30)(41,0)
\qbezier[15](22,0)(29.5,16)(37,0)
\qbezier[15](41,0)(49.5,16)(56,0)
\qbezier(12,0)(19.5,16)(27,0)
\qbezier(32,0)(39.5,16)(47,0)
\put(3,0){\vector(1,0){59}}
\end{picture} 
-
\begin{picture}(65,20)
\put(12,0){\circle*{2}}
\put(27,0){\circle*{2}}
\put(32,0){\circle*{2}}
\put(26,-4){\line(1,0){7}}
\put(47,0){\circle*{2}}
\qbezier[15](7,0)(14.5,16)(22,0)
\qbezier[20](3,0)(21,30)(41,0)
\qbezier[15](22,0)(29.5,16)(37,0)
\qbezier[15](41,0)(49.5,16)(56,0)
\qbezier(12,0)(22,16)(32,0)
\qbezier(27,0)(37,16)(47,0)
\put(3,0){\vector(1,0){59}}
\end{picture} 
\]

\medskip \noindent Notwithstanding, it remains worthwhile to investigate whether some weight system can be associated with distributions from noncommutative field theories.



\begin{thebibliography}{20}
\bibitem{quasiI} Bahns D, Doplicher S, Fredenhagen K and Piacitelli G 2005
 {\it Phys.\ Rev.} D {\bf 71} 025022
({\it Preprint} hep-th/0408204)
\bibitem{bahnsDiss} Bahns D 2004 {\it PhD Thesis} DESY-THESIS-2004-004
\bibitem{kassel} Kassel C 1995 {\it Quantum Groups} (New York: Springer)
\bibitem{barnatanIntro} Bar-Natan D 1995 {\it Topology} {\bf 34} 423
\bibitem{conneskreimer}
Connes A and Kreimer D 1999 {\it Eur. Phys. J.} C {\bf 7} 697
({\it Preprint} hep-th/9808042)
\bibitem{mmr} Bar-Natan D and Garoufalidis S 1996 {\it Inv. Math.} {\bf 125} 103 
\bibitem{stoimenow} Stoimenow A 1998 {\it PhD Thesis} FU Berlin 
\bibitem{quasiII}Bahns D, Doplicher S, Fredenhagen K and Piacitelli G, in preparation
\bibitem{bahnsIR} Bahns D, in preparation
\end{thebibliography}
\end{document}